\documentclass{article}
\usepackage{etoolbox}
\newtoggle{preprint}
\newtoggle{spl}
\newcommand{\preprint}[1]{\iftoggle{preprint}{#1}{}}
\newcommand{\spl}[1]{\iftoggle{spl}{#1}{}}

\togglefalse{spl}
\toggletrue{preprint}

\usepackage[hide=false,setmargin=true,marginparwidth=0.8in]{marginalia}
\newcommand{\manualendproof}{\hfill\qedsymbol\\[2mm]}

\spl{
\usepackage{amssymb}

\journal{Statistics and Probability Letters}

\input{global-macros}

\newcommand{\intdomain}{\mathcal{D}}

\newcommand{\N}{\mathbb{N}}

\newcommand{\R}{\mathbb{R}}

\newcommand{\absbig}[1]{\Big\lvert #1 \Big\rvert}
\newcommand{\abssmall}[1]{\lvert #1 \rvert}
\newcommand{\Norm}[1]{\left\lVert#1\right\rVert}

\newcommand{\modeaccent}[1]{\widehat{#1}}

\newcommand{\dist}{\pi} %

\newcommand{\llhood}{\ell_n}
\newcommand{\llhoodp}{\ell_n^{\paramidx}}
\newcommand{\remlep}{\hat \re_n^{\paramidx}}

\newcommand{\llandp}{{\ell^\dist_n}}

\newcommand{\logposthess}{\mb{\widehat H}_{n}}
\newcommand{\logposthessfunc}{\mb{H}_{n}}
\newcommand{\logposthessfuncp}{\mb{H}^{\paramidx}_{n}}

\newcommand{\approxdist}{\widetilde \dist} %
\newcommand{\LAapprox}{\approxdist_{\text{\upshape\tiny LA}}}

\newcommand{\dataidx}{\mb{Y}}

\newcommand{\data}{\mb{Y}^{(n)}} %
\newcommand{\paramtrue}{\mb{\theta}^{*}} %
\newcommand{\truelam}{\lambda^*}
\newcommand{\paramtrueidx}{\theta^{*}} %

\newcommand{\wtruetheta}{\omega^*_{\theta}}

\newcommand{\datatrueprob}{\PP^*}

\newcommand{\datatrueprobn}{\PP^*_{\!n}}

\newcommand{\poststar}[1]{\dist(#1 , \data)}

\newcommand{\parammode}{\modeaccent{\param}_n} %

\newcommand{\paramdim}{k}
\newcommand{\paramdimdum}{s}
\newcommand{\paramidxdim}{1}

\newcommand{\datadim}{d}

\newcommand{\param}{\mb{\theta}}

\newcommand{\paramidx}{\theta}

\newcommand{\paramspace}{\Theta}
\newcommand{\latentspace}{\Uu}
\newcommand{\latentvar}{\mb{U}}
\newcommand{\latentdim}{\paramdim}

\newcommand{\derivnum}{m} %
\newcommand{\derivbound}{M}

\newcommand{\dumderivvec}{\mb{\alpha}}
\newcommand{\dumderiv}{\alpha}

\newcommand{\llhoodmargin}{b}
\newcommand{\priorbig}{c_2}
\newcommand{\priorsmall}{c_1}

\newcommand{\eigen}{\lambda}

\newcommand{\conmargin}{\beta}

\newcommand{\hessbig}{\overline{\eta}}

\newcommand{\hesssmall}{\underline{\eta}}

\newcommand{\ball}[3]{B^{#1}_{#2}(#3)} %
\newcommand{\ballc}[3]{[B^{#1}_{#2}(#3)]^c} %

\newcommand{\constt}{C}

\newcommand{\universalradius}{\delta}
\newcommand{\universalradiusU}{\delta_\theta}

\newcommand{\mb}[1]{\boldsymbol{#1}}

\newcommand{\relerror}{\text{E}_{\mathrm{rel}}}

\newcommand{\iidsim}{\overset{iid}{\sim}}
\newcommand{\hess}{\mb{\widehat H}}

\newcommand{\ybar}{\overline{Y_n}}

\newcommand{\multiysum}[1]{N_{#1}}
\newcommand{\minmultiysum}{M}
\newcommand{\zbar}{\overline{Z_n}}

\newcommand{\re}{\omega}
\newcommand{\remle}{\hat\omega_{\theta}}

\newcommand{\multinomparamidx}{\psi}
\newcommand{\multidumparam}{\mb{\tau}}
\newcommand{\multidumparamidx}{\tau}

\newcommand{\multinomparam}{\mb{\multinomparamidx}}
\newcommand{\multinomorder}{k}

\newcommand{\simplex}{\mathcal{S}}
\newcommand{\ed}{\overset{d}{=}}

\newcommand{\modeparamidx}{\modeaccent{\paramidx}}
\newcommand{\modedumparamidx}{\modeaccent{\multidumparamidx}}
\newcommand{\modeparam}{\modeaccent{\param}}

\newtheorem{mytheorem}{Theorem}
\newtheorem{mylemma}{Lemma}
\newtheorem{mycorollary}{Corollary}
\newtheorem{myremark}{Remark}
\newtheorem{assumption}{Assumption}
\newtheorem{proposition}{Proposition}
}

\preprint{
\usepackage{
	amsmath,
	amsthm,
	amssymb,
	wrapfig,
	cases,
	mathtools,
	thmtools,
	array,
	bbm,
	bm,
	subfigure,
	makecell,
	esvect,
	mathrsfs,
	breqn,
	booktabs,
	upgreek,
	changepage,
	dsfont,
	fixme,
	listings,
	multirow,
	xargs,
	xstring,
	multicol,
	graphicx,
	float,
	color,
	enumerate,
	indentfirst,
	ifthen,
	wasysym,
	tikz,
	pgf,
	lmodern,
	titlesec,
	microtype
	}

\usepackage[utf8]{inputenc}
\usepackage[noend,ruled,vlined]{algorithm2e}
\usepackage{xpatch}
\makeatletter
\xpatchcmd{\algorithmic}
  {\ALG@tlm\z@}{\leftmargin\z@\ALG@tlm\z@}
  {}{}
\makeatother

\usepackage[hyperfootnotes=false, hidelinks]{hyperref}

\usepackage{natbib}
\setcitestyle{authoryear,round,citesep={;},aysep={,},yysep={;}}

\usetikzlibrary{shapes,decorations,arrows,calc,arrows.meta,fit,positioning}
\tikzset{
    -Latex,auto,node distance =1 cm and 1 cm,semithick,
    state/.style ={ellipse, draw, minimum width = 0.7 cm},
    point/.style = {circle, draw, inner sep=0.04cm,fill,node contents={}},
    bidirected/.style={Latex-Latex,dashed},
    el/.style = {inner sep=2pt, align=left, sloped}
}

\usepackage[capitalize,nameinlink]{cleveref}
\usepackage{crossreftools}
\pdfstringdefDisableCommands{%
    \let\Cref\crtCref
    \let\cref\crtcref
}

\hypersetup{%
    bookmarksnumbered, bookmarksopen=true, bookmarksopenlevel=1,%
}

\crefname{figure}{Figure}{Figures}
\crefname{subsection}{Subsection}{Subsections}
\crefname{lemma}{Lemma}{Lemmas}
\crefname{corollary}{Corollary}{Corollaries}
\crefname{theorem}{Theorem}{Theorems}
\crefname{assumption}{Assumption}{Assumptions}

\declaretheorem[name=Theorem]{theorem}

\declaretheorem[sibling=theorem,name=Lemma]{lemma}
\declaretheorem[sibling=theorem,name=Proposition]{proposition}
\declaretheorem[sibling=theorem,name=Corollary]{corollary}

\declaretheorem[name=Assumption]{assumption}
\declaretheorem[name=Assumption, numbered=no]{assumption*}

\declaretheorem[qed=$\triangleleft$,sibling=theorem,name=Remark]{remark}

\numberwithin{equation}{section}
\numberwithin{theorem}{section}
\numberwithin{lemma}{section}
\numberwithin{remark}{section}
\numberwithin{proposition}{section}
\numberwithin{definition}{section}
\numberwithin{corollary}{section}

\makeatletter
\renewcommand{\maketag@@@}[1]{\hbox{\m@th\normalsize\normalfont#1}}%
\makeatother

\makeatletter
\let\reftagform@=\tagform@
\def\tagform@#1{\maketag@@@{\ignorespaces\textcolor{gray}{(#1)}\unskip\@@italiccorr}}
\renewcommand{\eqref}[1]{\textup{\reftagform@{\ref{#1}}}}
\makeatother

\newcommand{\EE}{\mathbb{E}}

\newcommand{\NN}{\mathbb{N}}

\newcommand{\PP}{\mathbb{P}}

\newcommand{\RR}{\mathbb{R}}

\newcommand{\Uu}{\mathcal{U}}

\newcommand{\eps}{\varepsilon}

\def\[#1\]{\begin{equation}\begin{aligned}#1\end{aligned}\end{equation}}
\def\*[#1\]{\begin{equation*}\begin{aligned}#1\end{aligned}\end{equation*}}

\def\s*[#1\s]{\small\begin{align*}#1\end{align*}\normalsize}

\newcommand{\lcrx}[4][{-1}]{ 
	\IfEq{#1}{-1}{\left #2 {{{{#3}}}} \right #4}{
   	\IfEq{#1}{0}{#2 {{{{#3}}}} #4}{
	\IfEq{#1}{1}{\bigl #2 {{{{#3}}}} \bigr #4}{
	\IfEq{#1}{2}{\Bigl #2 {{{{#3}}}} \Bigr #4}{
	\IfEq{#1}{3}{\biggl #2 {{{{#3}}}} \biggr #4}{
	\IfEq{#1}{4}{\Biggl #2 {{{{#3}}}} \Biggr #4}{
    \GenericWarning{"4th argument to lcrx must be -1, 0, 1, 2, 3, or 4"}
    }}}}}}} %

\newcommand{\setdelim}{\ \vert \ }

\newcommand{\as}{\text{a.s.}}

\def\multiset#1#2{\ensuremath{\left(\kern-.3em\left(\genfrac{}{}{0pt}{}{#1}{#2}\right)\kern-.3em\right)}}

\DeclareMathOperator*{\argmin}{\arg\min} %
\DeclareMathOperator*{\argmax}{\arg\max} %
\DeclareMathOperator*{\newlim}{\mathrm{lim}\vphantom{\mathrm{infsup}}}
\DeclareMathOperator*{\newmin}{\mathrm{min}\vphantom{\mathrm{infsup}}}
\DeclareMathOperator*{\newmax}{\mathrm{max}\vphantom{\mathrm{infsup}}}
\DeclareMathOperator*{\newinf}{\mathrm{inf}\vphantom{\mathrm{infsup}}}
\DeclareMathOperator*{\newsup}{\mathrm{sup}\vphantom{\mathrm{infsup}}}
\renewcommand{\lim}{\newlim}
\renewcommand{\min}{\newmin}
\renewcommand{\max}{\newmax}
\renewcommand{\inf}{\newinf}
\renewcommand{\sup}{\newsup}

\newcommand{\Tr}{^{\scriptscriptstyle\text{T}}} %
\newcommand{\dee}{\mathrm{d}} %
\newcommand{\poissondist}{\mathrm{Pois}}

\newcommand{\bernoullidist}{\mathrm{Ber}}

\newcommand{\gammadist}{\mathrm{Gamma}}
\newcommand{\multidist}{\mathrm{MultiNom}}
\newcommand{\normaldist}{\mathrm{Gaussian}}

\newcommand{\uniformdist}{\mathrm{Unif}}

\newcommand{\sbra}[2][{-1}]{\lcrx[#1] [ {#2} ] }

\newcommand{\abs}[2][{-1}]{\lcrx[#1] \vert {#2} \vert }

\newcommand{\Nats}{\NN}

\newcommand{\Reals}{\RR}

\newcommand{\range}[2][{1}]{
	\IfEq{#1}{1}{\sbra{#2}}{\sbra{#2}_{#1}}}
\newcommand{\rangeO}[2][{0}]{
	\IfEq{#1}{0}{\sbra{#2}_0}{\sbra{#2}_{#1}}}

\title{On the Tightness of the Laplace Approximation\\ for Statistical Inference}
\author{Blair Bilodeau\\ University of Toronto \and Yanbo Tang\\ Imperial College London \and Alex Stringer\\ University of Waterloo}
\date{}

\setlength{\parindent}{0pt}
\setlength{\parskip}{6pt}
\usepackage{titlesec}
\newcommand{\addperiod}[1]{#1.}
\titleformat{\paragraph}[runin]
  {\normalfont\normalsize\bfseries}
  {\theparagraph}
  {1em}
  {\addperiod}
\titlespacing{\section}{0pt}{\parskip}{0pt}
\titlespacing{\subsection}{0pt}{\parskip}{0pt}
\titlespacing{\subsubsection}{0pt}{\parskip}{0pt}
\usepackage[letterpaper, margin=1in]{geometry}
}

\begin{document}

\spl{
\begin{frontmatter}

\title{On the Tightness of the Laplace Approximation\\ for Statistical Inference}

\author[inst1]{Blair Bilodeau}

\affiliation[inst1]{organization={Department of Statistical Sciences, University of Toronto},%
            addressline={27 King's College Circle}, 
            city={Toronto},
            postcode={M5S 1A4}, 
            state={Ontario},
            country={Canada}}

\author[inst2]{Yanbo Tang}
\author[inst3]{Alex Stringer}

\affiliation[inst2]{organization={Department of Mathematics, Imperial College London},%
            city={London SW7 2AZ},
            country={United Kingdom}}
            
\affiliation[inst3]{organization={Department of Statistics and Actuarial Science, University of Waterloo},%
            addressline={200 University Ave}, 
            city={Waterloo},
            postcode={N2L 3G1}, 
            state={Ontario},
            country={Canada}}

\begin{abstract}
Laplace's method is used to approximate intractable integrals in a statistical problems. The relative error rate of the approximation is not worse than $O_p(n^{-1})$. We provide the first statistical lower bounds showing that the $n^{-1}$ rate is tight.
\end{abstract}

\begin{highlights}
\item First stochastic lower bound on the error in the Laplace approximation in statistical integration problems.
\end{highlights}


\end{frontmatter}
}

\preprint{
\maketitle
\vspace{-25pt} %
\begin{abstract}
Laplace's method is used to approximate intractable integrals in a statistical problems. The relative error rate of the approximation is not worse than $O_p(n^{-1})$. We provide the first statistical lower bounds showing that the $n^{-1}$ rate is tight.
\end{abstract}
}

\section{Introduction}\label{sec:intro}

\subsection{Laplace Approximation in Statistical Problems}\label{subsec:laplaceapproximation}

Intractable integrals occur routinely in statistical modelling. Laplace's 
method 
\citep[e.g.,][Section~11.3.1]{davison_statistical_2008}
is a classical approach to approximating integrals. 
Two common statistical applications are approximating (a) the marginal likelihood in generalized mixed effects models
\citep{breslow_approximate_1993} 
and generalized additive models 
\citep{wood_fast_2011},
and (b) the normalizing constant required to compute posterior distributions
\citep{tierney_accurate_1986}
. 
It is known that the relative error in the Laplace approximation is no worse than $O_p(n^{-1})$ (where $O_p$ denotes \emph{stochastic boundedness}) under standard regularity conditions \citep{kass_validity_1990,bilodeau_stochastic_2021}. 
It is not known whether this $O_p(n^{-1})$ rate can be improved in general.
The contribution of this paper is to prove that the known $O_p(n^{-1})$ upper bound on the error is tight: to obtain a faster rate would require assumptions strong enough to rule out the most basic statistical models.

\subsection{Laplace's Method for Deterministic Integrals}\label{subsec:deterministic}

Given a fixed dimension $\paramdim \in \Nats$, domain $\intdomain \subseteq \Reals^\paramdim$, and function $h:\intdomain\to\Reals$, consider a sequence of (intractable) integrals $(I_n)_{n\in\Nats}$ defined by
$$
    I_{n} = \int_{\intdomain} e^{-nh(x)} \dee x.
$$
Laplace's method is to approximate $I_n$ by
\begin{equation}\label{eqn:laplacedeterministic}
    \tilde{I}_n = (2\pi)^{\paramdim/2}|\hess|^{-1/2}\exp\left\{-nh(\hat{u})\right\},
\end{equation}
where $\hat{u}=\argmin_{u\in\intdomain}h(u)$, and $\hess = \partial^{2}nh(\hat{u})$. 
Laplace first proposed the method in \citet{laplace1774memoire}.
\citet{olver68} first showed that the relative error $\relerror = \abssmall{I_n/\tilde{I}_n - 1}$ satisfies $\relerror \in O(n^{-1})$.
\citet{mcclure83} refined this analysis and extended it to two dimensions, while \citet{inglot_simple_2014} provided tighter constants and
the first lower bound showing $\relerror \in \Omega(n^{-1})$.

When Laplace's approximation is used in statistical problems, $n$ is taken to be a sample size, and $h$ depends on the negative log-likelihood function. 
The relative error $\relerror$ is then a random variable, 
and hence existing results for deterministic integrands do not apply in statistical problems. 

Obtaining bounds on the relative error
which hold uniformly over the randomness in the data is understood to be difficult. It is thus often preferred to seek stochastic convergence bounds which hold with high probability over the data.
For $\paramdim=1$, \citet{kass_validity_1990}
show that under various regularity conditions on the model the relative error 
is $O_p(n^{-1})$, where the stochasticity is with respect to the (unknown, potentially misspecified) data-generating measure.
A rigorous proof for $\paramdim\geq1$ can be found in \citet{bilodeau_stochastic_2021}.

\subsection{Contribution: Lower Bound in Statistical Problems}\label{subsec:contribution}
To the best of our knowledge, no lower bounds on the \emph{stochastic} relative error for the Laplace approximation are available.
A lower bound for all models is impossible: $\relerror=0$ almost surely for Gaussian models under any data-generating distribution as the Laplace approximation employs a quadratic approximation to the log-likelihood \citep[Section 11.3.1]{davison_statistical_2008}.
However, it remains to understand whether the $O_p(n^{-1})$ rate can be improved for any interesting class of models. The contribution of this paper is to answer this question in the negative: any set of assumptions under which a faster rate of convergence can be derived must be strong enough to rule out the most basic statistical models.

For deterministic integrals of the form discussed in \cref{subsec:deterministic}, \citet{inglot_simple_2014} show that under mild assumptions on $h$, there exists real constants 
$K_2,K_4, K_l, K_u$ 
and $\eta > 0$ such that
\begin{align*}
    \relerror = \abssmall{ I_n/\tilde{I}_n - 1} \geq \min\left\{ \abssmall{K_2/n + K_l/n^{1 + \eta/2} } , \abssmall{K_2/n + K_4/n^2 + K_u/n^{1 + \eta/2}} \right\},
\end{align*}
(taking $\lambda=n$ and $\alpha= 2$ in their notation). 
These constants depend on the integrand through $h$.
However, the existence of a lower bound on the relative error---depending on unknown, function-dependent constants---does not imply that the rate is tight in practice. There may be an interesting class of functions for which constants on lower-order terms in the expansion are zero, and hence a faster rate is attained. Further, this lower bound applies only to deterministic integrands, and not to statistical applications of Laplace's method.

In this paper, we take an alternative approach to proving that the $O_p(n^{-1})$ rate for the Laplace approximation is tight in statistical problems. Rather than proving a lower bound under assumptions that may be difficult to verify---leaving open the possibility of a faster rate being obtained in practice---we instead give three \emph{simple examples} of models in which the $O_p(n^{-1})$ rate is attained, including one in which the lower bound holds for any fixed parameter dimension $\paramdim\geq1$. We conclude that any set of assumptions under which a faster rate is obtained must exclude these models. We intentionally choose the simplest possible example models, to illustrate the strength of the assumptions that would be required in order to obtain a faster rate.

\section{Existing Result: Stochastic Upper Bound}\label{sec:assumptions}

\subsection{Laplace Approximation in Bayesian Inference}\label{subsec:upperbound}

Consider data $\data = (\dataidx_1,\dots,\dataidx_n) \subseteq \R^{\datadim}$ generated from some unknown joint probability measure $\datatrueprobn$. A parametric Bayesian model for $\data$
fixes a parameter space $\paramspace \subseteq \R^{\paramdim}$, a response density (likelihood) $\dist(\cdot \setdelim \param)$ indexed by $\param\in\paramspace$, and a prior density $\dist(\cdot)$ on $\paramspace$. Bayesian inferences are based on the posterior density, which itself depends on the \emph{normalizing constant}
\begin{equation}\label{eqn:normalizingconstant}
	\dist(\data) = \int_{\paramspace}\poststar{\param}\dee\param,
\end{equation}
where $\poststar{\param} = \dist(\data \setdelim \param)\dist(\param)$. The integral defining \cref{eqn:normalizingconstant} is intractable in general. In this paper, we focus on the use of the Laplace approximation to approximate \cref{eqn:normalizingconstant}. Similarly to \cref{subsec:deterministic}, define $\logposthessfunc(\param) = -\partial^{2}_{\param}\log\poststar{\param}$, $\parammode=\argmax{\poststar{\param}}$, and $\logposthess = \logposthessfunc(\parammode)$. The Laplace approximation to \cref{eqn:normalizingconstant} is
\begin{equation}\label{eqn:laplacebayesian}
\LAapprox(\data) = (2\pi)^{\paramdim/2}|\logposthess|^{-1/2}\poststar{\parammode}.
\end{equation}

A stochastic convergence rate for $\paramdim=1$ was derived by \citet{kass_validity_1990}. The $\paramdim\geq1$ case recently appeared as a corollary of \citet[Theorem 1]{bilodeau_stochastic_2021}, which we restate here:
\spl{\begin{mytheorem}}
\preprint{\begin{theorem}}[\citealt{bilodeau_stochastic_2021}]\label{thm:upperbound}
Suppose $\data\sim\datatrueprobn$ and that $\log\poststar{\param}$ is a $m\geq4$-times differentiable function of $\param$, $\datatrueprobn-a.s$. Under \cref{assn:kderiv,assn:hessian,assn:limsup,assn:consistency,assn:prior} (\cref{subsec:regularityconditions}), there exists $\constt>0$ with
$$
	\lim_{n\to\infty}\datatrueprobn\left( \absbig{\frac{\dist(\data)}{\LAapprox(\data)} - 1} \leq \constt \, n^{-1}\right) = 1.
$$
\spl{\end{mytheorem}}
\preprint{\end{theorem}}
\cref{thm:upperbound} provides a stochastic upper bound on the relative error in using $\LAapprox(\data)$ to approximate $\dist(\data)$: as long as the process $\datatrueprobn$ and model $\poststar{\param}$ together satisfy 
the regularity conditions outlined in \cref{subsec:regularityconditions} 
(which require neither a well-specified model nor i.i.d.\ observations), 
the relative error is $O_p(n^{-1})$.

\subsection{Laplace Approximation in Marginal Likelihood Models}\label{subsec:marginallikelihood}

Intractable integrals also appear regularly in marginal likelihood models. Once again, 
consider data $\data = (\dataidx_1,\dots,\dataidx_n) \subseteq \R^{\datadim}$ generated from some unknown joint probability measure $\datatrueprobn$. 
A marginal likelihood model for $\data$
fixes a parameter space $\paramspace \subseteq \R^{\paramdimdum}$ and a \emph{latent space} $\latentspace \subseteq \R^{\latentdim}$, a conditional response density $\dist(\cdot \setdelim \latentvar; \param)$ for a given $\param\in\paramspace$ and $\latentvar\in\latentspace$, and a \emph{latent generating density} $\dist(\cdot)$ on $\latentspace$. Inferences are based on the \emph{marginal likelihood}
\begin{equation}\label{eqn:marglikelihood}
	\dist(\data; \param) = \int_{\latentspace}\dist(\data \setdelim \latentvar; \param) \dist(\latentvar)\dee\latentvar.
\end{equation}

For any fixed $\param\in\paramspace$, this model is identical to the Bayesian model described in \cref{subsec:upperbound} by taking $\latentspace$ as the ``parameter space" to be integrated over and the latent-generating density as the ``prior". Under this relabelling, we can define the Laplace approximation exactly as in \cref{eqn:laplacebayesian}, and obtain the following.
\spl{\begin{mycorollary}}
\preprint{\begin{corollary}}\label{thm:upperbound-marg}
Fix $\param\in\paramspace$.
Suppose $\data\sim\datatrueprobn$ and that $\log\dist(\data,\latentvar;\param)$ is a $m\geq4$-times differentiable function of $\latentvar$, $\datatrueprobn-a.s$. Under \cref{assn:kderiv,assn:hessian,assn:limsup,assn:consistency,assn:prior} (\cref{subsec:regularityconditions}), there exists $\constt>0$ with
$$
	\lim_{n\to\infty}\datatrueprobn\left( \absbig{\frac{\dist(\data;\param)}{\LAapprox(\data;\param)} - 1} \leq \constt \, n^{-1}\right) = 1.
$$
\spl{\end{mycorollary}}
\preprint{\end{corollary}}

\section{Stochastic Lower Bounds}\label{sec:lower}

\subsection{Bayesian Inference: Coin Flips with Uniform Prior}\label{subsec:coinflipsregularity}

To establish the tightness of \cref{thm:upperbound}, we will prove a lower bound of $\Omega(n^{-1})$ almost surely for the following simple model:
\begin{equation}\label{eqn:bernoullimodel}
\begin{aligned}
Y_1,\ldots,Y_n \setdelim \paramidx \iidsim \bernoullidist(\paramidx), \ \paramidx\sim\uniformdist(0,1).
\end{aligned}
\end{equation}
That is, $\dist(\data \setdelim \paramidx) = \paramidx^{n\ybar}(1-\paramidx)^{n(1-\ybar)}$ where $\ybar = (1/n)\sum_{i=1}^{n}Y_i$, and $\dist(\paramidx) = I(0\leq\paramidx\leq1)$. 
We first verify the conditions of \cref{thm:upperbound} (\cref{assn:kderiv,assn:hessian,assn:limsup,assn:consistency,assn:prior} in \cref{subsec:regularityconditions}) hold when this model is well-specified (which only requires that the data are i.i.d., since binary data must have a Bernoulli distribution).
\begin{proposition}\label{prop:bernoulliassumptions}\label{PROP:BERNOULLIASSUMPTIONS}
For every $\paramtrueidx\in(0,1)$, \cref{eqn:bernoullimodel} satisfies \cref{assn:kderiv,assn:hessian,assn:limsup,assn:consistency,assn:prior} for $\datatrueprobn = \bernoullidist(\paramtrueidx)^{\otimes n}$.
\end{proposition}
The proof of \cref{prop:bernoulliassumptions} is given in \cref{sec:regularity-proof}. Having established that this simple model satisfies the assumptions under which \cref{thm:upperbound} holds, proving that it has a matching lower bound will establish the tightness of \cref{thm:upperbound}. 

\spl{\begin{mytheorem}}
\preprint{\begin{theorem}}
\label{thm:lowerbound}
Let $\paramtrueidx \in (0.25, 0.75)$ and $\datatrueprobn = \bernoullidist(\paramtrueidx)^{\otimes n}$.
For large enough $n$, \cref{eqn:bernoullimodel} satisfies
$$
    \absbig{\frac{\dist(\data)}{\LAapprox(\data)}-1}
    \geq \frac{1}{26n}
    \qquad \datatrueprobn-\as
$$
\spl{\end{mytheorem}}
\preprint{\end{theorem}}
\cref{thm:lowerbound} establishes that the $O_p(n^{-1})$ rate for the Laplace approximation 
is achieved by the simple Bernoulli model with Uniform prior. \cref{prop:bernoulliassumptions} states that this model satisfies \cref{assn:kderiv,assn:hessian,assn:limsup,assn:consistency,assn:prior}. Thus, any set of conditions assumed on $\datatrueprobn$ and $\poststar{\paramidx}$ must be strong enough to exclude the model of \cref{eqn:bernoullimodel} if a faster rate is to be obtained by the Laplace approximation. 

\subsection{Bayesian Inference: Multinomial with Dirichlet Prior}\label{subsec:multinomialregularity}

We demonstrate that the $n^{-1}$ rate is tight in a model with arbitrary fixed parameter dimension.
Consider the following standard multinomial model with uniform prior.
Fix $\multinomorder \in \Nats$, let $\simplex = \{\multinomparam\in[0,1]^\multinomorder: \sum_{j=1}^\multinomorder \multinomparamidx_j = 1\}$, and define
\[\label{eqn:multinommodel}
    \dataidx_1,\ldots,\dataidx_n \setdelim \multinomparam \iidsim \multidist(\multinomparam), \ \multinomparam\sim\uniformdist(\simplex).
\]
\begin{proposition}\label{prop:multinomialassumptions}
For every $\multinomorder\in\Nats$, the model given by \cref{eqn:multinommodel} reparameterized with $\paramidx_j = \log\multinomparamidx_j/\multinomparamidx_\multinomorder$ for $j\in[\multinomorder-1]$ satisfies \cref{assn:kderiv,assn:hessian,assn:limsup,assn:consistency,assn:prior}.
\end{proposition}

\spl{\begin{mytheorem}}
\preprint{\begin{theorem}}
\label{thm:multinom-lowerbound}
For every $\multinomorder\in\Nats$ and $\multinomparam\in\simplex$, if $\datatrueprobn = \multidist(\multinomparam)^{\otimes n}$,
then for large enough $n$,
the model given by \cref{eqn:multinommodel} reparameterized with $\paramidx_j = \log\multinomparamidx_j/\multinomparamidx_\multinomorder$ for $j\in[\multinomorder-1]$ satisfies
$$
    \absbig{\frac{\dist(\data)}{\LAapprox(\data)}-1}
    \geq \frac{1}{5n}
    \qquad \datatrueprobn-\as
$$
\spl{\end{mytheorem}}
\preprint{\end{theorem}}

\spl{\begin{myremark}}
\preprint{\begin{remark}}
While the Laplace approximation is not invariant to reparameterization, the upper bound from \cref{thm:upperbound} holds for any parameterization for which the regularity conditions are satisfied. \cref{thm:multinom-lowerbound} with $\multinomorder=2$ shows that the lower bound holds for an alternative parameterization of the Bernoulli model, and hence implies that this result is not an artifact of the particular parameterization chosen.
\preprint{\end{remark}}
\spl{\end{myremark}}

\subsection{Marginal Likelihood: Poisson with Gamma Random Effects}\label{subsec:poissontheorems}

We consider one of the simplest hierarchical models,
\begin{equation}\label{eqn:poissonmodel}
\begin{aligned}
Y_1,\ldots,Y_n \setdelim \re \iidsim \poissondist(\re\paramidx), \ \re\sim\gammadist(1,1); \ \paramidx>0.
\end{aligned}
\end{equation}
The model of \cref{eqn:poissonmodel} prescribes a joint likelihood $\dist(\data,\re;\paramidx)$, and hence a marginal likelihood $\dist(\data;\paramidx)=\int\dist(\data,\re;\paramidx)\dee\re$.
A Laplace approximate marginal likelihood $\LAapprox(\data;\paramidx)$ is obtained by applying \cref{eqn:laplacebayesian} to $\dist(\data;\paramidx)$ for every $\paramidx>0$. \cref{prop:poissonassumptions} states that \cref{thm:upperbound} applies to $\LAapprox(\data;\paramidx)$, pointwise in $\paramidx$.

\begin{proposition}\label{prop:poissonassumptions}\label{PROP:POISSONASSUMPTIONS}
For every $\paramidx>0$ and $\truelam>0$, \cref{eqn:poissonmodel} satisfies \cref{assn:kderiv,assn:hessian,assn:limsup,assn:consistency,assn:prior} for $\datatrueprobn = \poissondist(\truelam)^{\otimes n}$.
\end{proposition}

\cref{thm:poissonbound} now establishes that the rate is tight for $\LAapprox(\data;\paramidx)$ in the simple model of \cref{eqn:poissonmodel}, with conclusions similar to those of \cref{thm:lowerbound}.
\spl{\begin{mytheorem}}
\preprint{\begin{theorem}}
\label{thm:poissonbound}\label{THM:POISSONBOUND}
For all $\truelam > 0$ and large enough $n$, if $\datatrueprobn = \poissondist(\truelam)^{\otimes n}$ then \cref{eqn:poissonmodel} satisfies
$$
    \inf_{\paramidx>0}\absbig{\frac{\dist(\data ; \paramidx)}{\LAapprox(\data ; \paramidx)}-1}
    \geq \frac{1}{26n\truelam}
    \qquad \datatrueprobn-\as
$$
\spl{\end{mytheorem}}
\preprint{\end{theorem}}

\section*{Acknowledgements}
BB acknowledges support from the Vector Institute. YT acknowledges support from an Ontario Graduate Scholarship and the Vector Institute, and this work was partially completed while he was a PhD student at the University of Toronto.

\preprint{\bibliographystyle{abbrvnat}
\bibliography{references}}

\spl{\bibliographystyle{chicago}
\bibliography{references}}

\appendix
\spl{
\newpage
\begin{center}
{\Large\bf Supplementary Material for:\\ ``On the Tightness of the Laplace Approximation\\ for Statistical Inference"}
\end{center}

\setcounter{section}{0}
\setcounter{page}{1}
}

\section{Regularity Conditions}\label{subsec:regularityconditions}

We inherit the standard notation of \citet{bilodeau_stochastic_2021}. 
Let $\llhood(\param) = \log\dist(\data\setdelim\param)$ and $\llandp(\param) = \llhood(\param) + \log\dist(\param)$.
For any $f:\Reals^\paramdim\to\Reals$ and $\dumderivvec \in \Nats^\paramdim$, let 
\begin{equation*}
\begin{aligned}
    \partial^{\dumderivvec} f(\param)
    = \frac{\partial^{\abs{\dumderivvec}}}{\partial \param_1^{\dumderiv_1}\cdots\partial \param_\paramdim^{\dumderiv_\paramdim}} f(\param).
\end{aligned}
\end{equation*}
For $\paramdim\in\Nats$, $\mb{x}\in\Reals^\paramdim$, and $\delta>0$, let $\ball{\paramdim}{\mb{x}}{\delta} = \{\mb{x'}\in\Reals^\paramdim: \ \Norm{\mb{x}-\mb{x'}}_2 \leq \delta\}$ and $\ballc{\paramdim}{\mb{x}}{\delta} = \Reals^\paramdim \setminus \ball{\paramdim}{\mb{x}}{\delta}$.
Finally, the Eigenvalues of a $\paramdim$-dimensional square matrix $H$ are ordered $\eigen_\paramdim(H) \leq \cdots \leq \eigen_1(H)$. 

For a given data-generating distribution $\datatrueprobn$, we say the following assumptions hold if there exists $\universalradius>0$ and $\paramtrue\in\paramspace$ such that:

\begin{assumption}\label{assn:kderiv}
There exists $\derivnum\geq 4,\derivbound>0$ such that for all $\dumderivvec \in \N^\paramdim$ with $0 \leq \abssmall{\dumderivvec} \leq \derivnum $,
$$
	\lim_{n \to \infty} \datatrueprobn\Big[\sup_{\param \in \ball{\paramdim}{\paramtrue}{\universalradius}}\absbig{\partial^{\dumderivvec}\llandp(\param)} < n \derivbound\Big]=1.
$$
\end{assumption}

\begin{assumption} \label{assn:hessian}
There exist $0 < \hesssmall \leq \hessbig < \infty$ such that 
$$
	\lim_{n \to \infty} \datatrueprobn\Big[n \hesssmall \leq \inf_{\param \in \ball{\paramdim}{\paramtrue}{\universalradius}}\eigen_{\paramdim}(\logposthessfunc(\param)) \leq \sup_{\param \in \ball{\paramdim}{\paramtrue}{\universalradius}} \eigen_{1}(\logposthessfunc(\param))  \leq n \hessbig\Big] = 1.
$$
\end{assumption}

\begin{assumption}\label{assn:limsup}
There exists $\llhoodmargin>0$ such that 
$$
	\lim_{n \to \infty} \datatrueprobn\Big[\sup_{\param \in \ballc{\paramdim}{\paramtrue}{\universalradius}} \llhood(\param) - \llhood(\paramtrue) \leq -n \llhoodmargin\Big] = 1.
$$
\end{assumption}

\begin{assumption}\label{assn:consistency}
For any $\conmargin > 0$ and function $G(n)$ such that $\lim_{n \rightarrow \infty} G(n) = \infty$,
$$
  \lim_{n \to \infty} \datatrueprobn \left[ \frac{\sqrt{n}}{G(n)} \Norm{\parammode - \paramtrue}_2 > \conmargin \right] = 0.
$$
\end{assumption}

\begin{assumption}\label{assn:prior}
There exist  $0 < \priorsmall < \priorbig < \infty$ such that
$$
   \priorsmall \leq \inf_{\param \in \ball{\paramdim}{\paramtrue}{\universalradius}} \dist(\param) \leq \sup_{\param \in \ball{\paramdim}{\paramtrue}{\universalradius}} \dist(\param) \leq \priorbig.
$$
\end{assumption}

\section{Proofs for the Coin Flips with Uniform Prior Example\\ (Bayesian Inference)}\label{sec:regularity-proof}

To ease notation, let $\zbar = 1-\ybar$.
We also rely on the following classical lemma.
\spl{\begin{mylemma}}
\preprint{\begin{lemma}}
[Eqs. (1) and (2), \citet{robbins_stirling_1955}]
\label{fact:robbins}
For every $q\in\Nats$,
\begin{equation*}
    \left(2\pi q\right)^{1/2}\left(\frac{q}{e}\right)^{q}\exp\left(\frac{1}{12q+1}\right) 
    < \Gamma(q+1) < \left(2\pi q\right)^{1/2}\left(\frac{q}{e}\right)^{q}\exp\left(\frac{1}{12q}\right).
\end{equation*}
\preprint{\end{lemma}}
\spl{\end{mylemma}}

We now prove the main results from \cref{subsec:coinflipsregularity}.

\subsection{Proof of \cref{prop:bernoulliassumptions}}
For all $\paramidx\in (0,1)$,
\begin{equation*}
\begin{aligned}
    \llandp(\paramidx)
    = n\ybar\log\paramidx + n\zbar\log(1-\paramidx),
\end{aligned}
\end{equation*}
and for all $\alpha \in \Nats$
\begin{equation*}
\begin{aligned}
    \partial^\alpha \llandp(\paramidx)
    = (-1)^{\alpha+1}\frac{n\ybar}{\paramidx^\alpha} - \frac{n\zbar}{(1-\paramidx)^\alpha}.
\end{aligned}
\end{equation*}
Thus, for any $\paramtrueidx \in (0,1)$, taking $\universalradius = \min\{\paramtrueidx, 1-\paramtrueidx\} / 2$ implies that if $\alpha \leq \derivnum$,
\begin{equation*}
\begin{aligned}
    \sup_{\paramidx \in \ball{\paramidxdim}{\paramtrueidx}{\universalradius}}\abs{\partial^\alpha \llandp(\paramidx)}
    \leq 2n 
    \Big(\frac{2}{\min\{\paramtrueidx, 1-\paramtrueidx\}}\Big)^\derivnum,
\end{aligned}
\end{equation*}
by noting that $|\ybar| \leq 1$ and $|\zbar| \leq 1$. 
That is, \cref{assn:kderiv} holds by taking $\derivnum=4$.

Substituting in $\alpha = 2$,
\begin{equation*}
\begin{aligned}
    \logposthessfunc(\paramidx)
    = \frac{n\ybar}{\paramidx^2} + \frac{n\zbar}{(1-\paramidx)^2}.
\end{aligned}
\end{equation*}
By the strong law of large numbers, $\ybar \to \paramtrueidx$ and $\zbar \to 1-\paramtrueidx$ almost surely, so for all $\eps>0$ the following both hold almost surely for sufficiently large $n$ and for all $\paramidx \in \ball{\paramidxdim}{\paramtrueidx}{\universalradius}$ (with $\universalradius = \min\{\paramtrueidx, 1-\paramtrueidx\} / 2$):
\begin{equation*}
\begin{aligned}
    \logposthess(\paramidx)
    \leq n \left\{ 4\frac{(\paramtrueidx+\eps)}{(\paramtrueidx)^2} + 4\frac{(1-\paramtrueidx+\eps)}{(1-\paramtrueidx)^2}\right\}
\end{aligned}
\end{equation*}
and
\begin{equation*}
\begin{aligned}
    \logposthess(\paramidx)
    \geq n\left\{ \frac{(\paramtrueidx-\eps)}{4(\paramtrueidx)^2} + \frac{(1-\paramtrueidx-\eps)}{4(1-\paramtrueidx)^2}\right\}.
\end{aligned}
\end{equation*}
Since $\paramtrueidx \in (0,1)$, \cref{assn:hessian} is satisfied for $\hessbig$ and $\hesssmall$ depending on $1/\paramtrueidx$ and $1/(1-\paramtrueidx)$ appropriately. 

Next, for all $\paramidx\in(0,1)$,
\begin{align*}
    \llhood(\paramidx) - \llhood(\paramtrueidx)
    &= n\ybar\log(\paramidx/\paramtrueidx) + n\zbar((1-\paramidx)/(1-\paramtrueidx)).
\end{align*}
Again by the strong law of large numbers, this implies that 
\begin{align*}
    \frac{\llhood(\paramidx) - \llhood(\paramtrueidx)}{n}
    \longrightarrow
    -\mathrm{KL}(\bernoullidist(\paramtrueidx) \ \Vert \ \bernoullidist(\paramidx)) \qquad a.s.
\end{align*}
By Pinsker's inequality (e.g., Lemma~15.2 of \citet{wainwright}),
\begin{align*}
    \sup_{\paramidx \in \ballc{\paramidxdim}{\paramtrueidx}{\universalradius/2}}
    -\mathrm{KL}(\bernoullidist(\paramtrueidx) \ \Vert \ \bernoullidist(\paramidx))
    \leq 
    \sup_{\paramidx \in \ballc{\paramidxdim}{\paramtrueidx}{\universalradius/2}}
    \frac{-\abs{\paramidx-\paramtrueidx}^2}{2}
    \leq
    -\universalradius^2/8,
\end{align*}
so \cref{assn:limsup} holds.

Since the posterior is maximized at $\ybar$, \cref{assn:consistency} follows from the central limit theorem.

Finally, \cref{assn:prior} holds trivially.
\manualendproof

\subsection{Proof of \cref{thm:lowerbound}}
Let $\datatrueprob = \bernoullidist(\paramtrueidx)$, so that $\datatrueprobn = (\datatrueprob)^{\otimes n}$.
By definition,
\begin{equation*}
\begin{aligned}
    \LAapprox(\data) &= \left(\frac{2\pi}{n}\right)^{1/2}\Big(\ybar\Big)^{n\ybar+1/2}  \Big(\zbar\Big)^{n\zbar+1/2}.
\end{aligned}
\end{equation*}
The exact normalizing constant is
\begin{equation*}
\begin{aligned}
    \dist(\data) 
    &= \int_0^1 \dist(\data \setdelim \paramidx) \dist(\paramidx) \dee \paramidx
    = \int_0^1 \paramidx^{n \ybar} (1-\paramidx)^{n\zbar} \dee \paramidx
    = \frac{\Gamma(n\ybar+1)\Gamma(n\zbar+1)}{\Gamma(n+2)}.
\end{aligned}
\end{equation*}

Since $n\ybar$ and $n\zbar$ are both elements of $\Nats$,
by \cref{fact:robbins}
\*[
&\hspace{-1em}\frac{\Gamma(n\ybar+1)\Gamma(n\zbar+1)}{\Gamma(n+2)}  \\
&< (2\pi)^{1/2}\frac{n^{n+1}}{(n+1)^{n+3/2}}\exp\left(\frac{1}{12n\ybar \ \zbar}-\frac{1}{12n+13} +1 \right)\ybar^{n\ybar+1/2}\zbar^{n\zbar+1/2},
\]
and hence
\begin{equation*}
\begin{aligned}
    \frac{\dist(\data)}{\LAapprox(\data)} 
    &< \left(\frac{n}{n+1}\right)^{n+3/2}\exp\left( \frac{1}{12n\ybar \ \zbar}-\frac{1}{12n+13} +1\right) \\
    &= \exp(1)\left(1+\frac{1}{n}\right)^{-n-3/2}\exp\left(\frac{1}{12n\ybar \ \zbar} -\frac{1}{12n+13}\right).
\end{aligned}
\end{equation*}

Next, for all $n \in \Nats$ it holds that $1/(n+1) \leq \log(1+1/n)$, so rearranging implies that
\begin{equation*}
\begin{aligned}
    \frac{\dist(\data)}{\LAapprox(\data)} 
    &\leq \exp\Big(1 - \frac{n}{n+1} - \frac{3}{2(n+1)} -\frac{1}{12n+13} + \frac{1}{12n\ybar \ \zbar}\Big) \\
    &= \exp\Big(\frac{1}{12n\ybar \ \zbar} - \frac{14n+15}{(12n+13)(2n+2)}\Big).
\end{aligned}
\end{equation*}

Finally, if $\ybar \in (0.25,0.75)$ and $n>8$, then 
\*[
    \ybar \ \zbar
    \geq 0.1875
    > \frac{2n+2}{12n},
\]
and hence if $n>13$,
\*[
    \frac{1}{12n\ybar \ \zbar} - \frac{14n+15}{(12n+13)(2n+2)}
    &\leq \frac{1}{2n+2} - \frac{14n+15}{(12n+13)(2n+2)} \\
    &= - \frac{1}{12n+13} \\
    &\leq -\frac{1}{13n} \\
    &\leq \frac{\frac{-1}{26n}}{1-\frac{1}{26n}}.
\]
Using that $\frac{x}{1+x} \leq \log(1+x)$ for all $x > -1$, this implies that for sufficiently large $n$, if $\ybar \in (0.25,0.75)$ then
\*[
    \frac{\dist(\data)}{\LAapprox(\data)} 
    \leq 1-\frac{1}{26n}.
\]

Now, for the first time in the proof, we use the properties of the data-generating mechanism. In particular, by the strong law of large numbers it holds that for all $\eps>0$, $ \datatrueprob[\lim_{n\to\infty} \ybar \in (\paramtrueidx-\eps, \paramtrueidx+\eps)] = 1$, and thus since $\paramtrueidx \in (0.25, 0.75)$,
\begin{equation*}
\begin{aligned}
    \datatrueprob\Big[ 
    \exists N_0 \, \forall n > N_0 \quad
    \frac{\dist(\data)}{\LAapprox(\data)} 
    &\leq 1-\frac{1}{26n}
    \Big]
    =1.
\end{aligned}
\end{equation*}
\manualendproof

\section{Proofs for the Multinomial with Dirichlet Prior Example\\ (Bayesian Inference)}\label{appendix:multinomial}

\subsection{Reparameterized Model}\label{sec:reparam}
We begin with calculations for the Laplace approximation of the reparameterized model.

The normalizing constant is
\*[
\dist(\data) &= \int_{\multinomparam\in\simplex}\dist(\data , \multinomparam)\dee\multinomparam, \\
&= \Gamma(\multinomorder)\int_{\multinomparam\in\simplex} \prod_{j=1}^\multinomorder  
\multinomparamidx_j^{\multiysum{j}}\dee\multinomparam \\
&= \Gamma(\multinomorder)\frac{\prod_{j=1}^\multinomorder  \Gamma(\multiysum{j} + 1)}{\Gamma(n+\multinomorder)}.
\]

Using the standard log-odds reparametrization ($\paramidx_j = \log\multinomparamidx_j/\multinomparamidx_\multinomorder$ for $j=1,\ldots,\multinomorder-1$), we obtain $\multinomparamidx_j = \multinomparamidx_\multinomorder e^{\paramidx_j}$ and $\multinomparamidx_\multinomorder = 1/(1+\sum_{j=1}^{\multinomorder-1}e^{\paramidx_j})$. Note that the parameter space for $\param$ is $\Reals^{\multinomorder-1}$. The log-likelihood is thus
\*[
    \dist(\data\setdelim\param) = \exp\left\{ \sum_{j=1}^{\multinomorder-1}\multiysum{j}\paramidx_j - n\log\left( 1 + \sum_{j=1}^{\multinomorder-1}e^{\paramidx_j}\right)\right\}.
\]
Next, if $\multinomparam\sim\text{Dir}(1,\ldots,1)$ then
$\multinomparam\ed(X_1/X,\ldots,X_\multinomorder/X)$ where $X_1,\ldots,X_\multinomorder\iidsim\text{Exp}(1)$ and $X = \sum_{j=1}^{\multinomorder}X_j$. We therefore have 
$$
\left(e^{\paramidx_1},\ldots,e^{\paramidx_{\multinomorder-1}}\right)\ed\left(\frac{X_1}{X_\multinomorder},\ldots,\frac{X_{\multinomorder-1}}{X_\multinomorder}\right),
$$
and hence for every $\mb{\alpha} \in \Reals^{\multinomorder-1}$
\*[
F_{\param}(\alpha_1,\ldots,\alpha_{\multinomorder-1}) &= P\left(\paramidx_1\leq\alpha_1,\ldots,\paramidx_{\multinomorder-1}\leq\alpha_{\multinomorder-1}\right) \\
&= P\left(
\frac{X_1}{X_\multinomorder}\leq 
e^{\alpha_1},\ldots,
\frac{X_{\multinomorder-1}}{X_\multinomorder}\leq e^{\alpha_{\multinomorder-1}}
\right) \\
&= \EE_{X_\multinomorder}P\left( X_1\leq X_\multinomorder e^{\alpha_1},\ldots,X_{\multinomorder-1}\leq X_{\multinomorder} e^{\alpha_{\multinomorder-1}}\setdelim X_{\multinomorder}\right) \\
&= \EE_{X_{\multinomorder}}\prod_{j=1}^{\multinomorder-1}\exp\left( -X_{\multinomorder}e^{\alpha_j}\right) \\
&= \EE_{X_{\multinomorder}}\exp\left\{ -X_{\multinomorder}\sum_{j=1}^{\multinomorder-1}\exp\left(\alpha_j\right)\right\} \\
&= \frac{1}{1 + \sum_{j=1}^{\multinomorder-1}\exp\left( \alpha_j\right)}.
\]
The prior density is then
\*[
    \dist_{\param}(\mb{\alpha}) 
    &= 
    \frac
    {\partial^{\multinomorder-1}}
    {\partial\alpha_1\cdots\partial\alpha_{\multinomorder-1}}
    F_{\param}(\alpha_1,\ldots,\alpha_{\multinomorder-1}) \\
    &=
    \frac
    {\Gamma(\multinomorder)\exp\left(\sum_{j=1}^{\multinomorder-1}\alpha_j\right)}
    {\left(1 + \sum_{j=1}^{\multinomorder-1} e^{\alpha_j} \right)^{\multinomorder}}.
\] 
Rewriting the notation,
$$
\dist(\param) = \frac{\Gamma(\multinomorder)\exp\left(\sum_{j=1}^{\multinomorder-1}\paramidx_j\right)}{\left(1 + \sum_{j=1}^{\multinomorder-1}e^{\paramidx_j} \right)^{\multinomorder}}.
$$
The log joint posterior is thus
$$
\llandp(\param) = \log\dist(\param,\data) = \log\Gamma(\multinomorder) + \sum_{j=1}^{\multinomorder-1}\left(\multiysum{j}+1\right)\paramidx_j - (n+\multinomorder)\log\left( 1 + \sum_{j=1}^{\multinomorder-1}e^{\paramidx_j}\right).
$$
The first partial derivatives are
$$
\frac{\partial\llandp(\param)}{\partial\paramidx_j} = \multiysum{j}+1 - (n+\multinomorder)\frac{e^{\paramidx_j}}{1 + \sum_{l=1}^{\multinomorder-1}e^{\paramidx_l}},
$$
the posterior mode is
\*[
    \modeparamidx_j &= \log\frac{\multiysum{j}+1}{\multiysum{\multinomorder}+1},
\]
and the second order mixed partials are
\*[
\frac{\partial^2\llandp(\param)}{\partial\paramidx_i\partial\paramidx_j} &= 
\begin{cases}
-(n+\multinomorder)\frac{e^{\paramidx_j}}{1 + \sum_{l=1}^{\multinomorder-1}e^{\paramidx_l}}\left(1-\frac{e^{\paramidx_j}}{1 + \sum_{l=1}^{\multinomorder-1}e^{\paramidx_l}}\right) & i=j \\
(n+\multinomorder)\frac{e^{\paramidx_i}e^{\paramidx_j}}{\left(1 + \sum_{l=1}^{\multinomorder-1}e^{\paramidx_l}\right)^2} & i\neq j.
\end{cases}
\]

For notational convenience, let $\multidumparam \in [0,1]^{\multinomorder-1}$ be defined by
\*[
    \multidumparamidx_j
    = \frac{e^{\paramidx_j}}{1+\sum_{l=1}^{\multinomorder-1} e^{\paramidx_l}}.
\]

Then,
\*[
    \logposthessfunc(\param) = -\partial^2\llandp(\param) = (n+\multinomorder)\left\{\text{diag}\left(\multidumparam\right) - \multidumparam\multidumparam\Tr\right\},
\]
with determinant
\*[
    \abs{\logposthessfunc(\param)} 
    &= (n+\multinomorder)^{\multinomorder-1}|\text{diag}\left(\multidumparam\right)|\left(1 - \multidumparam\Tr\text{diag}\left(\multidumparam\right)^{-1}\multidumparam \right) \\
    &= (n+\multinomorder)^{\multinomorder-1}\left(\prod_{j=1}^{\multinomorder-1}\multidumparamidx_j\right)\left(1 - \sum_{i=1}^{\multinomorder-1}\sum_{j=1}^{\multinomorder-1}\multidumparamidx_i\multidumparamidx_j\text{diag}\left(\multidumparam\right)^{-1}_{ij}\right) \\
    &= (n+\multinomorder)^{\multinomorder-1}\left(\prod_{j=1}^{\multinomorder-1}\multidumparamidx_j\right)\left(1 - \sum_{j=1}^{\multinomorder-1}\multidumparamidx_j\right) \\
    &= (n+\multinomorder)^{\multinomorder-1}\left(\prod_{j=1}^{\multinomorder-1}\multidumparamidx_j\right)\frac{1}{1+\sum_{l=1}^{\multinomorder-1} e^{\paramidx_l}}.
\]

\subsection{Proof of \cref{prop:multinomialassumptions}}

We first verify \cref{assn:kderiv}.
Consider the rescaled log-likelihood combined with the log-prior
\[
\frac{\llandp(\param)}{n} = \frac{\log\Gamma(\multinomorder) + \sum_{j=1}^{\multinomorder-1}\left(\multiysum{j}+1\right)\paramidx_j - (n+\multinomorder)\log\left( 1 + \sum_{j=1}^{\multinomorder-1}e^{\paramidx_j}\right)}{n}.\label{eqn:resclaed_lik}
\]
First consider the derivatives of:
\*[
\frac{(n+\multinomorder)}{n} \log\left( 1 + \sum_{j=1}^{\multinomorder-1}e^{\paramidx_j}\right).
\]
Both $\log(1 + x)$ and $\sum_{i = 1}^{\multinomorder - 1}\exp(x_i)$ are real analytic functions, thus their composition must also be real analytic and therefore infinitely differentiable and have a convergent power series representation. In this case the power series converges for all values of $\param \in \Reals^{\multinomorder - 1}$, so for any compact set $A$ centered around $\paramtrue$, derivatives of total order $1,2,3$ and $4$ can be uniformly bounded by a constant $M_{\paramtrue} (A)$. 
As for the other terms in \cref{eqn:resclaed_lik}, the constant term does not appear in any derivatives, while the linear term satisfies
\*[\frac{\partial}{\partial \paramidx_i}\frac{\sum_{j=1}^{\multinomorder-1}\left(\multiysum{j}+1\right)\paramidx_j}{n} = \frac{\multiysum{i}+1}{n} \leq 2,\]
and is $0$ for any higher order derivatives.
Thus:
\*[\sup_{\param \in A} \absbig{\partial^{\dumderivvec}\llandp(\param)} \leq 2M_{\paramtrue}(A) + 2, \]
for $n > \multinomorder$ as then $(n+\multinomorder)/(n) \leq 2$, for $|\dumderivvec | \leq 4$. Taking $A$ to be the closure of the ball in \cref{assn:kderiv} then implies \cref{assn:kderiv} is satisfied.

As shown in \cref{sec:reparam}, $\logposthessfunc(\param) = (n+\multinomorder)\left\{\text{diag}\left(\multidumparam\right) - \multidumparam\multidumparam\Tr\right\}$, so for any $j$
\*[
    \frac{1}{n} \eigen_{j}(\logposthessfunc(\param))
    = \frac{n+\multinomorder}{n} \eigen_j(\text{diag}\left(\multidumparam\right) - \multidumparam\multidumparam\Tr).
\]
Further,
\*[
    \abs{\text{diag}\left(\multidumparam\right) - \multidumparam\multidumparam\Tr} 
    &= \left(\prod_{j=1}^{\multinomorder-1}\multidumparamidx_j\right)\frac{1}{1+\sum_{l=1}^{\multinomorder-1} e^{\paramidx_l}}.
\]
This determinant is strictly non-zero and is bounded for all possible values of $ \param \in \Reals^{\multinomorder - 1}$, so the eigenvalues are upper and lower bounded in any open ball, and thus \cref{assn:hessian} is satisfied.

Next, note that
\*[
\frac{\llhood(\param) - \llhood(\paramtrue)}{n} \xrightarrow{a.s.} 
\frac{\sum_{j=1}^{\multinomorder-1} \exp(\paramtrueidx_j)(\paramidx_j - \paramtrueidx_j)}{1 + \sum_{j=1}^{\multinomorder-1}e^{\paramtrueidx_j}}  - \log\left( 1 + \sum_{j=1}^{\multinomorder-1}e^{\paramidx_j}\right) + \log\left( 1 + \sum_{j=1}^{\multinomorder-1}e^{\paramtrueidx_j}\right).
\]
After some algebra it can be seen that the strictly convex function
\*[ 
    F(x) = \log\left(1 + \sum_{j=1}^{\multinomorder-1}e^{x_j} \right)
\]
satisfies
\begin{align*}
    \frac{\llhood(\param) - \llhood(\paramtrue)}{n} \xrightarrow{a.s.}
    -F(\param) + F(\paramtrue) + \nabla F(\paramtrue)^\top (\param -\paramtrue)
    =
    -D_F(\param, \paramtrue),
\end{align*}
where $D_F$ denotes the Bregman divergence.
Since Bregman divergence is strictly positive and convex in the first argument, this function is strictly negative for any value outside of a ball of radius $\delta$, and thus \cref{assn:limsup} holds.

Next, using the original parameterization,
 \*[
     \modeparamidx_j - \paramtrueidx_j
     &= \log\frac{\multiysum{j}+1}{\multiysum{\multinomorder}+1} - \log\frac{\multinomparamidx^*_j}{\multinomparamidx^*_k} \\
     &= \left(\log\frac{\multiysum{j}+1}{n}  - \log(\multinomparamidx^*_j)\right)- \left(\log\frac{\multiysum{\multinomorder}+1}{n} - \log(\multinomparamidx^*_k)\right).
 \]
By the multivariate CLT, the delta method (applied to $\log x$), and Slutsky's lemma we have that the joint vector converges to
\*[
 \begin{bmatrix}
 n^{1/2}\left(\log\frac{\multiysum{j}+1}{n}  - \log\multinomparamidx^*_j\right)\\
 n^{1/2}\left(\log\frac{\multiysum{\multinomorder}+1}{n} - \log\multinomparamidx^*_\multinomorder\right)
 \end{bmatrix} \xrightarrow{D} N\left( \begin{bmatrix} 0 \\ 0 \end{bmatrix} , \begin{bmatrix}
     \frac{1- \multinomparamidx^*_j}{\multinomparamidx^*_j} &1\\ 1 & \frac{1- \multinomparamidx^*_\multinomorder}{\multinomparamidx^*_\multinomorder}
 \end{bmatrix}\right).
 \]
Using Wald's device (i.e., for all vectors $a$ of compatible dimension, $a^\top X_n \xrightarrow{D} a^\top X$ if $X_n \xrightarrow{D} X$) with $a = (1, -1)^\top$ gives
 \*[ 
    n^{1/2}( \modeparamidx_j - \paramtrueidx_j) \xrightarrow{D} N\left(0,  \frac{1- \multinomparam_j}{\multinomparam_j} + \frac{1- \multinomparam_\multinomorder}{\multinomparam_\multinomorder} \right ).
 \]
Since the rate of consistency is $n^{-1/2}$, this implies \cref{assn:consistency}.

Finally, we note that the prior
 \*[
 \dist(\param) = \frac{\Gamma(\multinomorder)\exp\left(\sum_{j=1}^{\multinomorder-1}\paramidx_j\right)}{\left(1 + \sum_{j=1}^{\multinomorder-1}e^{\paramidx_j} \right)^{\multinomorder}},
 \]
 is strictly positive and upper/lower bounded in any open ball that is a subset of $\Reals^{\multinomorder-1}$, as the function is continuous and never takes on the value $0$ for any $\param \in \Reals^{\multinomorder-1}$.
 That is, \cref{assn:prior} is trivially satisfied.
\manualendproof

\subsection{Proof of \cref{thm:multinom-lowerbound}}

Substituting $\modeparam$ into the definition of $\multidumparam$, it follows that 
\*[
    \modedumparamidx_j = \frac{\multiysum{j}+1}{n+\multinomorder},
\]
and thus the determinant at the mode is
\*[
    \abs{\logposthess} 
    = (n+\multinomorder)^{\multinomorder-1}\prod_{j=1}^{\multinomorder}\frac{\multiysum{j}+1}{n+\multinomorder} 
    = \frac{1}{(n+\multinomorder)}\prod_{j=1}^{\multinomorder}\left(\multiysum{j}+1\right).
\]

The Laplace approximation is thus
\*[
\LAapprox(\data) &= (2\pi)^{(\multinomorder-1)/2}|\logposthess|^{-1/2}\poststar{\modeparam}\\
&= 
(2\pi)^{(\multinomorder-1)/2}(n+\multinomorder)^{1/2}
\Gamma(\multinomorder)
\prod_{j=1}^{\multinomorder}\left[\left(\multiysum{j}+1\right)^{-1/2}\right]\\
&\quad\quad \times 
\exp\left\{\sum_{j=1}^{\multinomorder-1}\left(\multiysum{j}+1\right)
\log\left(\frac{\multiysum{j}+1}{\multiysum{\multinomorder}+1}\right) -
(n+\multinomorder)\log\left(\frac{n+\multinomorder}{\multiysum{\multinomorder}+1}
\right)\right\} \\
&= 
(2\pi)^{(\multinomorder-1)/2}(n+\multinomorder)^{1/2}
\Gamma(\multinomorder)
\prod_{j=1}^{\multinomorder}\left[\left(\multiysum{j}+1\right)^{-1/2}\right]\\
&\quad\quad \times
\exp\left\{\sum_{j=1}^{\multinomorder}\left(\multiysum{j}+1\right)
\log\left(\multiysum{j}+1\right) -
(n+\multinomorder)\log\left(n+\multinomorder
\right)\right\}.
\]

Substituting this gives
\begin{align*}
    \frac{\dist(\data)}{\LAapprox(\data)} 
    &= 
    \frac{(2\pi)^{(1-\multinomorder)/2}}{\Gamma(n+\multinomorder) (n+\multinomorder)^{1/2}}
    \prod_{j=1}^{\multinomorder} \left[\left(\multiysum{j}+1\right)^{1/2}\Gamma(\multiysum{j}+1)\right] \\
    &\quad\quad \times
    \exp\left\{-\sum_{j=1}^{\multinomorder}\left(\multiysum{j}+1\right)
    \log\left(\multiysum{j}+1\right) +
    (n+\multinomorder)\log\left(n+\multinomorder
    \right)\right\}.
\end{align*}

Again by \cref{fact:robbins},
\begin{align*}
    &\hspace{-1em}\prod_{j=1}^{\multinomorder} \left[\left(\multiysum{j}+1\right)^{1/2}\Gamma(\multiysum{j}+1)\right] \\
    &\geq 
    \prod_{j = 1}^{\multinomorder}\left[ (2\pi )^{1/2} \left( \frac{\multiysum{j}}{\multiysum{j}  +1}\right)^{\multiysum{j}+1/2}\left(\multiysum{j}+1\right)^{\multiysum{j}+1} e^{-\multiysum{j}} \exp\left(\frac{1}{12\multiysum{j}+1}\right)\right]\\
    &= (2\pi)^{\multinomorder/2} \prod_{j = 1}^{\multinomorder}\left[ \left(1 + \frac{1}{\multiysum{j}}\right)^{-(1/2+\multiysum{j})} \exp\left( -\multiysum{j} + (\multiysum{j}+1)\log(\multiysum{j}+1) + \frac{1}{12\multiysum{j}+1}\right)\right]\\
    &=  (2\pi)^{\multinomorder/2} e^{-n} \prod_{j = 1}^{\multinomorder}\left(1 + \frac{1}{\multiysum{j}}\right)^{-(1/2+\multiysum{j})} \exp\left( \sum_{j=1}^\multinomorder (\multiysum{j}+1)\log(\multiysum{j}+1) + \sum_{j=1}^\multinomorder\frac{1}{12\multiysum{j}+1}\right). \\
\end{align*}
Similarly,
\begin{align*}
    &\hspace{-1em}\frac{1}{(n+\multinomorder)^{1/2}\Gamma(n+\multinomorder)} \\ 
    &\geq (2\pi)^{-1/2}(n+\multinomorder)^{-(n+\multinomorder)}\left(1-\frac{1}{n+\multinomorder}\right)^{1/2 - (n+\multinomorder)}e^{n+\multinomorder-1}\exp\left( -\frac{1}{12(n+\multinomorder-1)}\right).
\end{align*}
Combining these bounds gives
\begin{align*}
    \frac{\dist(\data)}{\LAapprox(\data)} 
    &\geq 
    \left(1-\frac{1}{n+\multinomorder}\right)^{1/2-(n+\multinomorder)}
    \prod_{j = 1}^{\multinomorder}\left(1 + \frac{1}{\multiysum{j}}\right)^{-(1/2+\multiysum{j})} \\
    &\quad\quad \times
    \exp\left(\multinomorder-1+
    \sum_{j=1}^\multinomorder\frac{1}{12\multiysum{j}+1}
    - \frac{1}{12( n +\multinomorder - 1 )}
    \right).
\end{align*}

We now use that $\log(1-1/x) \leq -1/x$ and $\log(1+1/x)\leq 1/x$ to obtain
\*[
    \frac{\dist(\data)}{\LAapprox(\data)} 
    &\geq
    \exp\left(\frac{n+\multinomorder-1/2}{n+\multinomorder}\right) \\
    &\quad \times
    \prod_{j = 1}^{\multinomorder}\exp\left(-\frac{1/2+\multiysum{j}}{1+\multiysum{j}}\right)
    \exp\left(\multinomorder-1+
    \sum_{j=1}^\multinomorder\frac{1}{12\multiysum{j}+1}
    - \frac{1}{12( n +\multinomorder - 1 )}
    \right) \\
    &=
    \exp\left(\frac{1}{2(n+\multinomorder)}
    +\sum_{j=1}^\multinomorder\frac{1}{2(1+\multiysum{j})}
    + \sum_{j=1}^\multinomorder\frac{1}{12\multiysum{j}+1}
    - \frac{1}{12( n +\multinomorder - 1 )}
    \right) \\
    &\geq
    \exp\left(\frac{1}{2(n+\multinomorder)}
    +\frac{1}{2(1+\minmultiysum)}
    + \frac{1}{12\minmultiysum+1}
    - \frac{1}{12( n +\multinomorder - 1 )}
    \right),
\]
where $\minmultiysum = \min_{j\in[\multinomorder]} \multiysum{j}$. Note that 
\*[
\frac{1}{2(n+\multinomorder)}  - \frac{1}{12(n+\multinomorder - 1)} \geq \frac{1}{3(n+\multinomorder - 1)} - \frac{1}{12(n+\multinomorder - 1)} = \frac{1}{4(n+\multinomorder - 1)}
\]
since $(n + k - 1)/(n+k ) \geq 2/3$ for $\multinomorder \geq 2$ and $n \geq 1$.
Combining this with the fact that the other terms involved in the exponent are positive,
\*[\frac{\dist(\data)}{\LAapprox(\data)} 
    &\geq \exp\left( \frac{1}{4(n+\multinomorder - 1)}\right) \geq 1 + \frac{1}{4(n+\multinomorder - 1)}, 
\]
as $\exp(x) \geq 1 + x$ for $x \geq 0$. 
Thus, for $n > 4(\multinomorder-1)$, we obtain
\*[
    \frac{\dist(\data)}{\LAapprox(\data)}
    \geq 1 + \frac{1}{5n}.
\]
\manualendproof

\section{Proofs for the Poisson with Gamma Random Effects Example\\ (Marginal Likelihood)}\label{sec:poissonproof}

\subsection{Random Effects Model}\label{sec:remodel}
The joint likelihood is
$$
    \dist(\data,\re;\paramidx) = \frac{\left(\re\paramidx\right)^{n\ybar}\exp\left\{-\re\left(n\paramidx+1\right)\right\}}{\prod_{i=1}^{n}\Gamma\left(Y_i+1\right)}.
$$
Thus, it is possible to compute the exact marginal likelihood
$$
    \dist(\data;\paramidx) = \int\dist(\data,\re;\paramidx)\dee\re = \frac{\Gamma(n\ybar+1)}{\prod_{i=1}^{n}\Gamma\left(Y_i+1\right)}\frac{\paramidx^{n\ybar}}{\left(n\paramidx+1\right)^{n\ybar+1}}.
$$
For each $\paramidx > 0$, let $\llhoodp(\re) = \log\dist(\data,\re;\paramidx)$, so that
\begin{equation}\label{eqn:llhoodp}
\begin{aligned}
    \llhoodp(\re)
    &= n\ybar \log(\re\paramidx) - \re(n\paramidx+1) - \sum_{i=1}^n \log(Y_i!), \\
    \frac{\dee}{\dee \re} \llhoodp(\re)
    &= \frac{n\ybar}{\re} - (n\paramidx+1), \\
    \frac{\dee^2}{\dee^2 \re} \llhoodp(\re)
    &= -\frac{n\ybar}{\re^2}.
\end{aligned}
\end{equation}
Thus, $\remlep = n\ybar / (n\paramidx+1)$, and so
$$
    \abs{\logposthessfuncp(\remlep)}^{-1/2}
    = \frac{\sqrt{n\ybar}}{n\paramidx+1}.
$$

\subsection{Proof of \cref{PROP:POISSONASSUMPTIONS}}
Fix a value of $\paramidx > 0$.
In \cref{assn:kderiv,assn:hessian,assn:limsup,assn:consistency,assn:prior} replace $\paramtrue$ with $\wtruetheta = \truelam/\paramidx$, $\universalradius$ with $\universalradiusU = \truelam/(2\paramidx)$ and $\parammode$ with $\remlep$.

Using the expressions from \cref{eqn:llhoodp}, for $\re \in (\wtruetheta-\universalradiusU, \wtruetheta+\universalradiusU)$ and $k>1$, the strong law of large numbers implies
\begin{align*}
    \frac{1}{n}\left|\frac{\dee}{\dee \re} \llhoodp(\re)\right|
    &=\left| \frac{\ybar}{\re} - (\paramidx+1/n)\right| \xrightarrow{a.s.} \left| \frac{\truelam}{\re} - \paramidx\right| \leq 
    \paramidx, \\
    \frac{1}{n}\left|\frac{\dee^j}{\dee^j \re} \llhoodp(\re)\right|
    &= \frac{(j - 1)\ybar}{\re^j} \xrightarrow{a.s.} \frac{(j - 1)\truelam}{\re^j} \leq \frac{(j-1)(2\paramidx)^j}{(\truelam)^{j-1}}.
\end{align*}
Thus, \cref{assn:kderiv} holds for a constant appropriately depending on $\paramidx$ and $\truelam$.
Similarly, for $\re \in (\wtruetheta-\universalradiusU, \re+\universalradiusU)$
\begin{align*}
    -\frac{1}{n}\frac{\dee^2}{\dee^2 \re} \llhoodp(\re)
    &= \frac{\ybar}{\re^2} \xrightarrow{a.s.} \frac{\truelam}{\re^2} \in \frac{1}{\truelam}\Big((4/9)\paramidx^2, 4\paramidx^2\Big),
\end{align*}
so \cref{assn:hessian} is satisfied.

For Assumption 3, again by the strong law of large numbers,
\begin{align*}
    \frac{\llhoodp(\re) - \log\dist(\re) - \llhoodp(\wtruetheta) + \log\dist(\wtruetheta)}{n} \xrightarrow{a.s.} \log\left(\frac{\re}{\wtruetheta}\right) + \frac{(\wtruetheta - \re)}{\wtruetheta} = -D_{F}(\re; \wtruetheta),
\end{align*}
where $D_{F}(\re, \wtruetheta)$ is the Bregman divergence induced by the function $F(x) = -\log(x)$ on the convex set $(0, \infty)$. \cref{assn:limsup} 
then follows from the non-negativity of the Bregman divergence and its strict convexity in the first argument.

Solving the score equation to compute the maxima gives
\begin{equation*}
\begin{aligned}
    \sqrt{n}\Big(\remle - \wtruetheta\Big)
    &= \sqrt{n}\left(\frac{\ybar}{\paramidx+1/n} - \frac{\truelam}{\paramidx}\right)
    = \frac{\sqrt{n}}{\paramidx+1/n} \cdot \left((\ybar - \truelam) - \frac{\truelam}{n\paramidx} \right)
    \rightsquigarrow \normaldist(0, \paramidx^{-2}),
\end{aligned}
\end{equation*}
where the last step follows by the central limit theorem.
Thus, \cref{assn:consistency} holds, and finally
\cref{assn:prior} is satisfied by the strict positivity of the gamma density on any compact subset of $(0,\infty)$.
\manualendproof

\subsection{Proof of \cref{thm:poissonbound}}
This implies that the Laplace approximate marginal likelihood is
\begin{equation*}\begin{aligned}
    \LAapprox(\data;\paramidx) &=
    (2\pi)^{1/2} \cdot \abs{\logposthessfuncp(\remlep)}^{-1/2} \cdot
    \dist(\data,\re; \paramidx), \\
    &= \frac{\left(2\pi\right)^{1/2}}{\prod_{i=1}^{n}\Gamma(Y_i+1)}\left(n\ybar\right)^{n\ybar+1/2}\frac{\paramidx^{n\ybar}}{\left(n\paramidx+1\right)^{n\ybar+1}}\exp\left(-n\ybar\right),
\end{aligned}\end{equation*}
so for each $\paramidx>0$,
$$
    \frac{\dist(\data;\paramidx)}{\LAapprox(\data;\paramidx)} = \frac{\Gamma(n\ybar+1)}{\left(2\pi\right)^{1/2}\big(n\ybar\big)^{n\ybar+1/2} \cdot\exp\left(-n\ybar\right)}.
$$
Note that this ratio does not depend on $\paramidx$, yielding the uniform statement in the theorem. 
Again applying 
\cref{fact:robbins},
we have
\begin{equation*}
\begin{aligned}
    \frac{\dist(\data;\paramidx)}{\LAapprox(\data;\paramidx)} &> 
    \frac{\left(2\pi\right)^{1/2}\big(n\ybar\big)^{n\ybar+1/2}\exp\left(-n\ybar + \frac{1}{12n\ybar+1}\right)}
    {\left(2\pi\right)^{1/2}\big(n\ybar\big)^{n\ybar+1/2} \cdot\exp\left(-n\ybar\right)}, \\
    &= \exp\left(\frac{1}{12n\ybar+1}\right) \\
    &\geq 1 + \frac{1}{12n\ybar+1},
\end{aligned}
\end{equation*}
where the last step uses that $e^x \geq 1+x$ for all $x$.

As in the proof of \cref{thm:lowerbound}, we now use a property of $\datatrueprob = \poissondist(\truelam)$. In particular, if $n \ybar > 1$ and $\ybar \leq 2\truelam$, then
\*[
    \frac{1}{12n\ybar+1}
    \geq \frac{1}{13n\ybar}
    \geq \frac{1}{26n\truelam}.
\]
By the strong law of large numbers, $\ybar \to \truelam$ almost surely (and hence $n \ybar \to \infty$ almost surely), which completes the proof.
\manualendproof

\end{document}